\def\year{2019}\relax
\documentclass[12pt]{article} 
\usepackage{times}  
\usepackage{helvet}  
\usepackage{courier}  
\usepackage{url}  
\usepackage{graphicx}  
\usepackage{appendix}
\usepackage{amsmath,fullpage, graphicx, enumerate, verbatim}
\usepackage{natbib}
\usepackage{amssymb}
\usepackage{algorithm}
\usepackage{dsfont}
\usepackage{algorithmic}
\usepackage{blindtext}
\usepackage{booktabs}  
\usepackage{threeparttable}
\usepackage{multicol}  
\usepackage{multirow}
\newtheorem{assumption}{\bf Assumption}[section]
\newtheorem{lemma}{\bf Lemma}[section]
\newtheorem{corollary}{\bf Corollary}[section]
\newtheorem{thm}{\bf Theorem}[section]
\newtheorem{defn}{\bf Definition}[section]

\title{Convergence Rates of Posterior Distributions in Markov Decision Process}
\author{Zhen Li, Eric Laber \\Department of Statistics\\ NC State University}
\date{}
 
\begin{document}
%

\maketitle
\begin{abstract}
In this paper, we show the convergence rates of posterior distributions of the model dynamics in a MDP for both episodic and continuous tasks. The theoretical results hold for general state and action space and the parameter space of the dynamics can be infinite dimensional. Moreover, we show the convergence rates of posterior distributions of the mean accumulative reward under a fixed or the optimal policy and of the regret bound. A variant of Thompson sampling algorithm is proposed which provides both posterior convergence rates for the dynamics and the regret-type bound. Then the previous results are extended to Markov games. Finally, we show numerical results with three simulation scenarios and conclude with discussions.
\end{abstract}

\section{Introduction}
\label{sec:intro}
A Markov Decision Process (MDP) is a discrete time stochastic control process which provides tools to model sequential decision making process under uncertainty \citep{howard1960dynamic}. MDPs have been widely applied in the fields of economies, precision medicine, robotics and games \citep{sutton1998introduction}\citep{alagoz2010markov}\citep{littman1994markov}. If the dynamics (environment) of a MDP is known, we can determine the optimal policy to obtain the most desirable mean outcome in the whole process with the method of reinforcement learning \citep{sutton1998introduction}. Moreover, the inference of the dynamics improve interpretability of the selected policy which is of great interest in the fields of mobile health and precision medicine. Thus a problem of interest in MDPs is to quantify the uncertainty of its dynamics (e.g. the transition probability and reward function). \cite{abbasi2011improved} and \cite{bastani2017mostly} show the confidence sets for the arm parameters under greedy algorithms in contextual bandits. \cite{jaksch2010near} gives the confidence intervals for the transition probability and reward function under a greedy algorithm in MDPs. \cite{theocharous2017posterior} shows the convergence rate of the estimator of parameters indexing the transition probability under a variant of Thompson sampling algorithm. However, most of the previous methods require that the state space $\mathcal{S}$ and action space $\mathcal{A}$ are finite and the parameter space of dynamics $\Theta$ is finite dimensional. In this paper, we apply a Bayesian approach to show the convergence rates of posterior distributions of the dynamics for general measurable space $\mathcal{S}$, $\mathcal{A}$ and infinite dimensional space $\Theta$ under any policy. Moreover, we show the convergence rates of posterior distributions of the mean accumulative reward under a fixed or the optimal policy and of the regret bound.

Thompson sampling is an online decision making algorithm which balances the exploitation and exploration \citep{thompson1933likelihood}. In the algorithm, a learner samples the dynamics from the posterior distribution at each time point given the past observations, compute the optimal policy based on the sampled dynamics and take action under the policy. In contrast to frequentist optimism algorithms which use the point estimator for the dynamics \citep{lai1985asymptotically}\citep{bartlett2009regal}, Thompson sampling has many advantages regarding both performance and computational issue \citep{osband2017posterior}. Several Thompson sampling algorithms for MDPs have been proposed including TSMDP \citep{gopalan2015thompson}, TSDE \citep{ouyang2017learning}, DS-PSRL \citep{theocharous2017posterior} and the regret bound for these algorithms are given. However, they need either finiteness of $\mathcal{S}$, $\mathcal{A}$ or restriction of the dimension of $\Theta$ and most of them do not provide the inference of the model dynamics. In this paper, we propose a Thompson algorithm which provides both convergence rates of posterior distributions for the dynamics and the regret bound for the number of suboptimal action selected.

In Section 2, we introduce the notation and present assumptions and theorems for posterior convergence rates in single-agent MDPs. In Section 3, we propose a variant of Thompson sampling algorithm and give a theorem regarding the convergence rate and regret bound of the algorithm. In Section 4, we extend the previous results to Markov games. In Section 5, we present numerical results of our proposed method with three simulation scenarios. Finally, we conclude this work with discussions and comments in Section 6.

\section{Posterior Convergence Rates in Markov Decision Process}
We consider a Markov Decision Process (MDP) $(\mathcal{S},\mathcal{A},r,p_\theta)$, where $\mathcal{S}$ is the state space, $\mathcal{A}$ is the action space and $r:\mathcal{S}\times \mathcal{A}\rightarrow \mathbb{R}$ is the reward function. $p_\theta:\mathcal{S}^2\times \mathcal{A}\rightarrow \mathbb{R}$ denotes the transition density such that $p_\theta(S_{t+1}|S_t,A_t)$ is the conditional density of $S_{t+1}\in \mathcal{S}$ given $S_t\in \mathcal{S}$ and $A_t\in \mathcal{A}$, $t\in \mathcal{T}=\{0,1,2,\ldots\}$, $\theta\in \Theta$. The parameter space $\Theta$ can be infinite dimensional and the true parameter of the transition density is denoted as $\theta_0$. We define the history up to time $t$ as $H_t=(S_0,A_0,R_0,S_1,A_1,R_1\ldots,S_{t-1},A_{t-1},R_{t-1},S_t)$ where $R_t=r(S_t,A_t)$ is the reward and the action $A_t$ depends on $H_t$, $t\in \mathcal{T}$.  Let $\mathcal{B}_{\mathcal{A}}$ denote the set of
distributions over $\mathcal{A}$. A policy, $\pi$, is an infinite sequence of functions $\pi_t:\mathrm{dom}\; H_t\rightarrow \mathcal{B}_\mathcal{A}$ such that under the policy $\pi$, a decision maker presented with the history $H_t$ at time $t$ will select an action $A_t$ from the distribution of $\pi_t(H_t)$, with density $p_{\pi_t}(A_t|H_t)$. We define
$V_\theta(\pi)=\mathbb{E}_{\theta,\pi}(\sum_{t=0}^\infty \gamma^t R_t)$ where $\mathbb{E}_{\theta,\pi}$ denotes the expectation w.r.t. the distribution induced by the dynamics $\theta$ and the policy $\pi\in \Pi$ and $\gamma \in (0,1]$ is a discount factor.
Given a class of policies $\Pi$ and $\theta\in \Theta$, the optimal policy $\pi_\theta^*\in \Pi$ satisfies $V_\theta(\pi_\theta^*)\geq V_\theta(\pi)$ for all $\pi\in \Pi$.

Suppose under a policy $\pi$ we observe $H_t$ , the complete history up to time $t$. We give a prior $P_t$ to $\theta\in \Theta$ and the posterior distribution is denoted as $P_t(\cdot|H_t)$ and we study the convergence rate of $P_t(\cdot|H_t)$. Define a semi-distance on $\Theta$: \begin{equation}
    \label{dist}
    d_\mu^2(\theta,\theta')=\int h^2(p_\theta(\cdot|s,a),p_{\theta'}(\cdot|s,a))\mu(s,a)d\lambda(s,a),
\end{equation}
where $h^2(p_1,p_2)=\int (\sqrt{p_1(x)}-\sqrt{p_2(x)})^2 dm(x)$ is the Hellinger distance between $p_1(\cdot)$ and $p_2(\cdot)$ w.r.t. a measure $m$ on $\mathcal{S}$; $\mu$ is a non-negative function and $\lambda$ is a reference finite measure on $\mathcal{S}\times \mathcal{A}$. We use $\lfloor x\rfloor$ to denote the greatest integer less than or equal to $x$. Let $P_{\theta,\pi}$ denote the probability distribution induced by the dynamics $\theta$ and the policy $\pi$. Define $N(\epsilon,\Theta,d_\mu)$ as the $\epsilon$-covering number which is the minimal number of balls with radius $\epsilon$ to cover $\Theta$ with respect to $d_\mu$. Moreover, let $K(p_1; p_2)$ be the Kullback-Leibler divergence so that $K(p_1;p_2)=\mathbb{E}_{p_1}\log \frac{p_1(X)}{p_2(X)}$ and define a discrepancy measure by $V_n(p_1;p_2)=\mathbb{E}_{p_1} (\log \frac{p_1(X)}{p_2(X)}-K(p_1;p_2))^n$. Then we let $$
\begin{aligned}
B_t(\theta_0,\epsilon;n)=&\{\theta\in \Theta: K(p_{\theta_0,\pi}^{(t)};p_{\theta,\pi}^{(t)})\leq t\epsilon^2,\\
&V_n(p_{\theta_0,\pi}^{(t)};p_{\theta,\pi}^{(t)})\leq t^{n/2}\epsilon^n\}
\end{aligned}
$$ where $p^{(t)}_{\theta,\pi}$ denotes the density of $H_t$ under the dynamics $\theta$ and the policy $\pi$.

To obtain the posterior convergence rate in a MDP, we need to first construct a test of $\theta_0$ versus $\{\theta:d_\mu(\theta,\theta_0)>\epsilon\}$, the complement of the ball around $\theta_0$. The idea of constructing a test for proof of the convergence rate comes from \cite{le2012asymptotic}, \cite{lecam1973convergence}, \cite{ghosal2000convergence}, \cite{ghosal2007convergence} and many tests were proposed for different problems but little is under the MDP scenario. Such a test can be found given Lemma \ref{testError} and some entropy conditions for the space $\Theta$. Before giving Lemma \ref{testError}, we first give an assumption for the function $\mu$ in (\ref{dist}). Define $\widetilde{H}_t=(H_t,A_t,R_t)$ and denote the conditional probability distribution of $(S_{i+j},A_{i+j})$ given $\widetilde{H}_i$ under the dynamics $\theta$ and the policy $\pi$ as $P^{\widetilde{H}_i,j}_{\theta,\pi}(\cdot,\cdot)$, with density $p_{\theta,\pi}^{\widetilde{H}_i,j}$, $i,j\geq 0$, $\theta\in \Theta$. Then Assumption \ref{A4} is defined as follows.
    \begin{assumption}
\label{A4}
For a policy $\pi\in \Pi$, there exists two positive integers $k$ and $l$, a non-increasing sequence of positive numbers $\{\alpha_i\}_{i\geq 0}$ and two non-negative functions $\mu$, $\nu$ such that for $\forall \theta\in \Theta$, $\forall i\geq 0$ and all possible $\widetilde{H}_i$, 
\begin{equation}
\label{mu}
     \alpha_i \mu(\cdot,\cdot)\leq  \frac{1}{k}\sum_{j=1}^k p^{\widetilde{H}_i,j}_{\theta,\pi}(\cdot,\cdot),
\end{equation}
\begin{equation}
\label{nu}
    p^{\widetilde{H}_i,l}_{\theta,\pi}(\cdot,\cdot)\leq \nu(\cdot,\cdot).
\end{equation}
\end{assumption}

Assumption \ref{A4} shows a kind of uniformity for the MDP. In the case where $\mathcal{S}$ and $\mathcal{A}$ is discrete, (\ref{mu}) implies every pair $(s,a)\in \mathcal{S}\times \mathcal{A}$ should be recurrent if $\mu$ is a positive function. That is reasonable since we hope to visit any pair $(s,a)\in \mathcal{S}\times\mathcal{A}$ infinitely often to learn the transition density $p_\theta$ fully.

If there exist $d>c>0$ such that $c\leq p^{\widetilde{H}_i,1}_{\theta,\pi}(\cdot,\cdot)\leq d$ for $\forall \widetilde{H}_i$ and $\forall \theta\in \Theta$, then it is obvious to see that Assumption \ref{A4} holds with $k=l=1$, $\alpha_i=1$ for $\forall i\geq 0$, $\mu=c$ and $\nu=d$.
Next we give a weaker assumption under which Assumption \ref{A4} holds in the case where $|\mathcal{S}|$ and $|\mathcal{A}|$ is finite, which is stated in Theorem \ref{weakA4}.
\begin{defn}
Suppose $|\mathcal{S}|<\infty$. For $\forall s,s'\in \mathcal{S}$, we say $(s,s')$ is connected w.r.t. $P_{\theta,\pi}$ if for $\forall i\geq 0$, $\exists\; j\geq 1$ such that $P_{\theta,\pi}(S_{i+j}=s'|S_i=s,\tilde{H}_{i-1})>0$ for $\forall \tilde{H}_{i-1}$.
\end{defn}
\begin{assumption}
\label{A5}
There exist $d>c>0$ such that $c\cdot p_{\theta_0}(s'|s,a)\leq p_\theta(s'|s,a)\leq d\cdot p_{\theta_0}(s'|s,a)$ for $\forall s,s'\in \mathcal{S}$, $\forall a\in \mathcal{A}$ and $\forall \theta\in \Theta$.
\end{assumption}
    \begin{thm}
    \label{weakA4}
    Suppose $|\mathcal{S}|<\infty$, $|\mathcal{A}|<\infty$ and Assumption \ref{A5} is satisfied. If (i) $\exists\; \delta_t>0$ such that $p_{\pi_t}(a|H_t)\geq \delta_t$ for $\forall a\in \mathcal{A}$ and $\forall H_t$ and (ii) $(s,s')$ is connected w.r.t. $P_{\theta,\pi}$ for $\forall s,s'\in \mathcal{S}$, $\forall \theta\in \Theta$, then there exist positive functions $\mu$ and $\nu$ satisfying Assumption \ref{A4}.
    \end{thm}
    
In Theorem \ref{weakA4}, (i) is satisfied if $\epsilon$-greedy algorithm is applied, allowing for decaying $\epsilon$, and (ii) guarantees every state in $\mathcal{S}$ is visited infinitely often which is necessary to learn the 
transition density $p_\theta$ fully. Assumption \ref{A5} implies the transition density $p_\theta(\cdot|s,a)$ possesses the same support for all $\theta\in \Theta$ and should be uniformly bounded for all $\theta\in \Theta$. This is equivalent to the bound condition for the log-likelihood ratio given in \cite{gopalan2015thompson}. The compactness of $\Theta$ and the smoothness of $p_\theta$ can be a sufficient condition for the existence of the bound and we will demonstrate in more detail later.

With Assumption \ref{A4}, we give Lemma \ref{testError} which helps constructing a test of $\theta_0$ and the complement of the ball around $\theta_0$.
    \begin{lemma}
\label{testError}
Under Assumption \ref{A4}, there exist constants $K>0$, $\xi\in (0,1)$ and tests $\phi_t$ such that
\begin{equation}
\label{typeI}
\mathbb{E}_{\theta_0,\pi}\phi_t \leq \exp\{-K\sum_{i=0}^{\lfloor t/(k+l)\rfloor-1} \alpha_{i(k+l)} d_\mu^2(\theta_0,\theta_1)\},
\end{equation}
\begin{equation}
\label{typeII}
\underset{B}{\sup}\;\mathbb{E}_{\theta,\pi}(1-\phi_t)\leq \exp\{-K\sum_{i=0}^{\lfloor t/(k+l)\rfloor-1} \alpha_{i(k+l)} d_\mu^2(\theta_0,\theta_1)\},
\end{equation}
no matter what the distribution of $X_0$ is and $B=\{\theta\in\Theta:d_\nu(\theta,\theta_1)\leq \xi \alpha_{(\lfloor t/(k+l)\rfloor-1)(k+l)}^{1/2} d_\mu(\theta_0,\theta_1)\}$.
\end{lemma}
Applying Lemma \ref{testError}, we can obtain the theorem of posterior convergence rates in MDP.
\begin{thm}
\label{MDP}
Under Assumption \ref{A4}, let $\epsilon_t>0$, $\epsilon_t\rightarrow 0$, $\bar{\epsilon}_t=(\sum_{i=0}^{\lfloor t/(k+l)\rfloor-1} \alpha_{i(k+l)}\epsilon_t^2/t)^{\frac{1}{2}}$, $1/(t\bar{\epsilon}_t^2)=\mathcal{O}(1)$ such that for every sufficient large $j\in \mathbb{N}$,
\begin{equation}
\label{entropyUB}
\begin{aligned}
    \underset{\epsilon>\epsilon_t}{\sup}\log &N(\frac{1}{2}\alpha_{(\lfloor t/(k+l)\rfloor-1)(k+l)}^{1/2}\xi \epsilon,\\
    &\{\theta\in \Theta:d_\mu(\theta,\theta_0)<\epsilon\},d_\nu)\leq t\bar{\epsilon}_t^2;
    \end{aligned}
\end{equation}
\begin{equation}
\label{priorLB}
    \frac{P_t(\theta\in \Theta:j\epsilon_t<d_\mu(\theta,\theta_0)\leq (j+1)\epsilon_t)}{P_t(B_t(\theta_0,\bar{\epsilon}_t;n))}\leq \exp\{Kt\bar{\epsilon}_t^2 j^2/2\}.
    \end{equation}
Then for every $M_t\rightarrow \infty$, we have that
\begin{equation}
\label{conv}
    P_t(\theta\in \Theta:d_\mu(\theta,\theta_0)> M_t\epsilon_t|H_t)\rightarrow 0
\end{equation}
in $P_{\theta_0,\pi}$-probability, where $\xi$ and $K$ are the same as those in Lemma \ref{testError}.
\end{thm}

Here $\epsilon_t$ is known as the posterior convergence rate in \cite{ghosal2017fundamentals}. The condition (\ref{entropyUB}) bounds the complexity of the parameter space $\Theta$ and it is combined with Lemma \ref{testError} to construct a test of $\theta_0$ versus $\{\theta:d_\mu(\theta,\theta_0)>\epsilon\}$. The condition (\ref{priorLB}) guarantees that the prior $P_t$ possesses enough mass around the true parameter $\theta_0$, i.e. $B_t(\theta_0,\bar{\epsilon}_t;n)$. Notice that the theorem holds no matter what the policy $\pi$ is.

Furthermore, suppose we can compute the V-function of a policy, $V_\theta(\pi)$, given the dynamics $\theta$ is known, applying methods such as planning and dynamic programming \citep{sutton1998introduction}. Also, we can compute the optimal policy $\pi_\theta^*$ under the dynamics $\theta$ applying methods such as Q-learning and value iteration \citep{sutton1998introduction}. Then $V_\theta(\pi)$ and $V_\theta(\pi_\theta^*)$ can be reasonable estimators for the true V-functions $V_{\theta_0}(\pi)$ and $V_{\theta_0}(\pi_{\theta_0}^*)$ respectively. $\pi_\theta^*$ is a reasonable estimator for the true optimal policy $\pi_{\theta_0}^*$ and we are interested in the regret bound of $V_{\theta_0}(\pi_{\theta_0}^*)-V_{\theta_0}(\pi_{\theta}^*)$. Under some mild assumptions, we can obtain the posterior convergence rates for the estimated V-functions and the regret bound in Corollary \ref{regret}.
\begin{assumption}
\label{A6}
$\lambda(\mathcal{S},\mathcal{A})<\infty$ and there exist positive numbers $C_1$, $C_2$, $C_3$, $C_4$, such that $\gamma C_1C_2\lambda(\mathcal{S},\mathcal{A})<1$ and
\begin{enumerate}[(i)]
    \item $p_{\pi_t}(a|H_t)\leq C_1$ for $\forall a\in \mathcal{A}$, $\forall \pi\in \Pi$ and $\forall H_t$, $t\in \mathcal{T}$;
    \item $p_{\theta}(s'|s,a)\leq C_2$ for $\forall s',s\in \mathcal{S},\forall a\in \mathcal{A}, \forall \theta\in \Theta$;
    \item $|p_\theta(s'|s,a)-p_{\theta_0}(s'|s,a)|\leq C_3 d_\mu(\theta,\theta_0)^\rho$ for $\forall s',s\in \mathcal{S}, a\in \mathcal{A},\theta\in \Theta$ and some $\rho\in (0,1]$;
    \item $|r(s,a)|\leq C_4$ for $\forall s\in \mathcal{S},\forall a\in \mathcal{A}$.
\end{enumerate}
\end{assumption}
\begin{corollary}
\label{regret}
Suppose Assumption \ref{A6} and (\ref{conv}) hold, then for $\forall \tilde{\pi}\in \Pi$ and every $M_t\rightarrow \infty$, we have
\begin{equation}
\label{Vpi}
P_t(\theta\in \Theta: |V_{\theta}(\tilde{\pi})-V_{\theta_0}(\tilde{\pi})|> M_t \epsilon_t^\rho|H_t)\rightarrow 0,
\end{equation}
\begin{equation}
    \label{Vpiopt1}
    P_t(\theta\in \Theta:|V_{\theta}(\pi_\theta^*)-V_{\theta_0}(\pi_{\theta_0}^*)|> M_t \epsilon_t^\rho|H_t)\rightarrow 0,
\end{equation}
\begin{equation}
\label{Vpiopt}
P_t(\theta\in \Theta:|V_{\theta_0}(\pi_\theta^*)-V_{\theta_0}(\pi_{\theta_0}^*)|> M_t \epsilon_t^\rho|H_t)\rightarrow 0
\end{equation}
in $P_{\theta_0,\pi}$-probability.
\end{corollary}

In Assumption \ref{A6}, (i), (ii) and (iv) are satisfied if $|\mathcal{S}|,|\mathcal{A}|<\infty$. (iii) states that the difference between the densities corresponding to $\theta$ and $\theta_0$ should be controlled by the distance $d_\mu(\theta,\theta_0)$ in some sense. Then we can obtain the posterior convergence rate $\epsilon_t^\rho$ for the estimated V-functions and the regret bound.

If the MDP can be broken into episodes \citep{sutton1998introduction} and we consider a class of policies $\Pi$ in which $\pi_t(H_t)$ depends only on the history of the current episode, Assumption \ref{A4} is no longer necessary and ideal functions $\mu$ and $\nu$ can be found to measure the distance between different elements in $\Theta$ (\ref{dist}). In a episodic task, each episode starts from a state $S_0$ sampled from a certain distribution and ends in a terminal state at time $T$, which can be random. Suppose $N$ episodes are completed at time $t$ and the posterior convergence rate is now in terms of $N$ instead of the overall time points $t$. As before, we will first give a lemma to construct a test with errors decreasing exponentially with $N$.

\begin{lemma}
\label{testError2}
Suppose $N$ episodes are completed at time $t$ and Assumption \ref{A5} holds. Then there exists constants $K>0$ and $\xi\in (0,1)$, functions $\mu$ and $\nu$, and tests $\phi_t$ such that
$$
\mathbb{E}_{\theta_0,\pi}\phi_t \leq \exp\{-KN d_\mu^2(\theta_0,\theta_1)\},
$$
$$
\underset{\theta\in\Theta:d_\nu(\theta,\theta_1)\leq \xi d_\mu(\theta_0,\theta_1)}{\sup}\mathbb{E}_{\theta,\pi}(1-\phi_t)\leq \exp\{-KN d_\mu^2(\theta_0,\theta_1)\},
$$
no matter what the distribution of $X_0$ is.
\end{lemma}
The functions $\mu$ and $\nu$ in Lemma \ref{testError2} both have the form $\sum_{t=0}^\infty \delta^t p_{\theta_0,\pi}^{t}(\cdot,\cdot)$ with different $\delta\in (0,1)$, where $p_{\theta_0,\pi}^{t}(\cdot,\cdot)$ is the density of $P_{\theta_0,\pi}(S_t\in \cdot,A_t\in \cdot)$ w.r.t. $\lambda$. The measure $\sum_{t=0}^\infty \delta^t P_{\theta_0,\pi}(\cdot,\cdot)$ on a set $S\times A\in \mathcal{S}\times\mathcal{A}$ is positive whenever $P_{\theta_0,\pi}(S_t\in S,A_t\in A)>0$ for some $t$. Moreover, the measure is a weighted sum of probability measures with different time points and early time points exert more influence on the measure. That is to say, the more likely the state and action falls into $S$ and $A$ in an early time point, the larger the value of the measure on this set is. Define
$$
\begin{aligned}
&B^*(\theta_0,\epsilon;n)=\{\theta\in \Theta: \mathbb{E}_{\theta_0,\pi} \sum_{t=1}^T \log \frac{p_{\theta_0}(S_t|S_{t-1},A_{t-1})}{p_{\theta}(S_t|S_{t-1},A_{t-1})}\leq \epsilon^2,\\
&\mathbb{E}_{\theta_0,\pi} |\sum_{t=1}^T \log \frac{p_{\theta_0}(S_t|S_{t-1},A_{t-1})}{p_{\theta}(S_t|S_{t-1},A_{t-1})}|^n\leq \epsilon^n\}.
\end{aligned}
$$
Then we can obtain Theorem \ref{MDP2} for posterior convergence rates in an episodic task.
\begin{thm}
\label{MDP2}
Suppose $N$ episodes are completed at time $t$ and Assumption \ref{A5} holds. Let $\epsilon_t>0$, $\epsilon_t\rightarrow 0$, $1/(N\epsilon_t^2)=\mathcal{O}(1)$ such that for every sufficient large $j\in \mathbb{N}$,
\begin{equation}
\label{entropyUB2}
    \underset{\epsilon>\epsilon_t}{\sup}\log N(\frac{1}{2}\xi \epsilon,\{\theta\in \Theta:d_\mu(\theta,\theta_0)<\epsilon\},d_\nu)\leq N\epsilon_t^2;
\end{equation}
\begin{equation}
\label{priorLB2}
    \frac{P_t(\theta\in \Theta:j\epsilon_t<d_\mu(\theta,\theta_0)\leq (j+1)\epsilon_t)}{P_t(B^*(\theta_0,\epsilon_t;n))}\leq \exp\{KN\epsilon_t^2 j^2/2\}.
    \end{equation}
Then for every $M_t\rightarrow \infty$, we have that
\begin{equation}
    P_t(\theta\in \Theta:d_\mu(\theta,\theta_0)\geq M_t\epsilon_t|H_t)\rightarrow 0
\end{equation}
in $P_{\theta_0,\pi}$-probability, where $\xi$, $K$, $\mu$ and $\nu$ are the same as those in Lemma \ref{testError2}.  
\end{thm}


\section{$\epsilon$-greedy Thompson Sampling}
Thompson sampling is an online decision-making algorithm balancing the trade-off between the exploration and exploitation \citep{thompson1933likelihood}. Based on Theorem \ref{MDP}, we propose an algorithm named $\epsilon$-greedy Thompson sampling which ensures both the concentration of the dynamics $\theta$ and the regret-type bound during the process. At each time point $t$, a prior for $\theta$ is given and we compute the posterior distribution of $\theta$, i.e. $P_t(\cdot|H_t)$. Then we sample $\theta$ from $P_t(\cdot|H_t)$, compute the optimal policy $\pi_\theta^*$ and select the action applying the $\epsilon$-greedy method. The $\epsilon$-greedy method is employed for every action to be sampled infinite times so the transition density $p_\theta$ can be learned fully. The whole process of $\epsilon$-greedy Thompson sampling is given in Algorithm \ref{TS}.
\begin{algorithm}[H]
	\begin{algorithmic}[1]
	\STATE Suppose we have a sequence of non-negative numbers $\{\delta_t\}_{t\in\mathcal{T}}$ and a sequence of priors $\{P_t\}_{t\in\mathcal{T}}$.
		\FOR {iteration $t=0,1,2,...$} 
		\STATE Sample $u\sim U(0,1)$.
		\STATE If $u>\delta_t$, (i) sample $\theta \sim P_t(\cdot|H_t)$; (ii) compute the optimal policy $\pi_\theta^*$ and (iii) select $A_t\sim \pi_{\theta,t}^*(H_t)$.\\
		Otherwise select $A_t$ from all the actions in $\mathcal{A}$ with equal probability.
		\STATE Observe $R_t,S_{t+1}$.
		\ENDFOR
	\end{algorithmic}
	\caption{$\epsilon$-greedy Thompson sampling}
	\label{TS}
\end{algorithm}
To obtain the regret bound of Algorithm \ref{TS}, we consider a class of \textit{stationary} policies $\Pi$ in which $\pi_t(H_t)$ only depends on the current state $S_t$. Then it is known that $\pi_{\theta,t}^*$ does not depend on $t$ where $\pi_{\theta}^*$ is the optimal policy in this class under $\theta$ \citep{sutton1998introduction}. Define the Q-function 
$$
Q_{\theta}(s,a,\pi)=\mathbb{E}_{\theta,\pi}[\sum_{t=0}^\infty \gamma^t R_t|S_0=s,A_0=a]
$$
and we know $\pi_{\theta,t}^*(H_t)=\underset{a\in\mathcal{A}}{\text{argmax}}Q_{\theta}(S_t,a,\pi_\theta^*)$. Then under some conditions, we can ensure both the concentration of the dynamics $\theta$ to $\theta_0$ and the regret-type bound for Algorithm \ref{TS} in Theorem \ref{Thompson}.
\begin{thm}
\label{Thompson}
For the $\epsilon$-greedy Thompson sampling process given in Algorithm \ref{TS}, suppose the assumptions in Theorem \ref{MDP} and Assumption \ref{A6} hold and $\Pi$ is a class of stationary policies. Moreover, suppose we have (i) $\bar{\epsilon}_t\geq t^{-\beta}$ for some $\beta\in (0,\frac{1}{2})$ and $n(1-2\beta)>2$ in Theorem \ref{MDP}; (ii) $\exists \;b>0$ such that $Q_{\theta_0}(s,\pi_{\theta_0,t}^*(s),\pi_{\theta_0}^*)\geq Q_{\theta_0}(s,a,\pi_{\theta_0}^*)+b$ for $\forall s\in \mathcal{S}$ and $\forall a\neq \pi_{\theta_0,t}^*(s)$; (iii) there exists a non-decreasing function $f$ such that $\sum_{i=0}^t \delta_i\leq f(t)$ and $f(t)\rightarrow \infty$ as $n\rightarrow \infty$.  Then for $\forall \delta\in (0,1)$, there exists $T_0\in \mathbb{N}$ such that for $\forall T\geq T_0$ and $\forall d>0$, we have
\begin{equation}
    \label{regretBound}
P_{\theta_0}(\sum_{t=0}^T 1(A_t\neq \pi_{\theta_0,t}^*(S_t))\leq f(T)^{1+d}+B(\delta))\geq 1-\delta,
\end{equation}
where $B(\delta)$ does not depend on $T$ and for sufficiently large $M$,
\begin{equation}
\label{convTS}
    P_t(\theta\in \Theta:d_\mu(\theta,\theta_0)> M\epsilon_T|H_T)\rightarrow 0
\end{equation}
almost surely as $T\rightarrow \infty$, where $\epsilon_T$ is given in Theorem \ref{MDP}.
\end{thm}
The condition (i) and (ii) guarantee the optimal policy $\pi_\theta^*$ we get at each time $t$ equals $\pi_{\theta_0}^*$, with sufficiently large probability and (ii) is satisfied when $|\mathcal{S}|,|\mathcal{A}|<\infty$. $(\ref{regretBound})$ gives a regret bound for the number of time points in $1,\ldots,T$ when the optimal action is not selected. The regret bound is determined by $f(\cdot)$ which depends on the exploration rate $\delta_t$, $t\in \mathcal{T}$. However, the exploration rate $\delta_t$ should be large enough to guarantee the transition density $p_\theta(s'|s,a)$ can be learned fully so the convergence result (\ref{convTS}) holds. Thus Algorithm \ref{TS} can be seen as a trade-off between selecting the optimal action and learning the dynamics of the environment.

\section{Posterior Convergence Rates in Markov Games}
In Section 2, we consider the Markov decision process with a single agent. A MDP can be generalized to the multi-agent scenario which is known as the stochastic game. A stochastic game with $K$ agents is represented by a tuple $(\mathcal{S},\mathcal{A}_1,\ldots,\mathcal{A}_K,r_1,\ldots,r_K,p_\theta)$, where $\mathcal{S}$ is the state space, $\mathcal{A}_k$ is the action space of the $k$th agent and $r_k:\mathcal{S}\times \mathcal{A}_1\times \ldots\times \mathcal{A}_K\rightarrow \mathbb{R}$ is the reward function of the $k$th agent, $k=1,\ldots,K$. $p_\theta:\mathcal{S}^2\times\mathcal{A}_1\times\ldots\times \mathcal{A}_K\rightarrow \mathbb{R}$ denotes the transition density such that $p_\theta(S_{t+1}|S_t,A_{1t},\ldots,A_{Kt})$ is the conditional density of $S_{t+1}\in \mathcal{S}$ given $S_t\in \mathcal{S}$ and $A_{kt}\in \mathcal{A}_k$, $k=1,\ldots,K$. The parameter of the transition density $\theta$ belongs to the space $\Theta$ which is of our interest and $\theta_0$ denotes the true parameter. Denote $\boldsymbol{A}_t=(A_{1t},\ldots,A_{Kt})$ as the joint action of the $K$ agents at time $t$ and $\boldsymbol{R}_t=(R_{1t},\ldots,R_{Kt})$ where $R_{kt}=r_k(S_t,A_{1t},\ldots,A_{Kt})$ is the reward of the $k$th agent at time $t$, $t\in \mathcal{T}$. Then the history up to time $t$ is defined as $H_t=(S_0,\boldsymbol{A}_0,\boldsymbol{R}_0,\ldots,S_{t-1},\boldsymbol{A}_{t-1},\boldsymbol{R}_{t-1},S_t)$. The policy of the $k$th agent, $\pi_k$ is an infinite sequence of functions $\pi_{kt}:\mathrm{dom}\; H_t\rightarrow \mathcal{B}_{\mathcal{A}_k}$ and the $k$th agent presented with the history $H_t$ at time $t$ will select an action $A_{kt}$ from the distribution of $\pi_{kt}(H_t)$, $k=1,\ldots,K$. For $k=1,\ldots,K$, we define $V_{k\theta}(\pi_1,\ldots,\pi_K)=\mathbb{E}_{\theta,\pi_1,\ldots,\pi_K}(\sum_{t=0}^\infty \gamma^t R_{kt})$ where $\mathbb{E}_{\theta,\pi_1,\ldots,\pi_K}$ denotes the expectation w.r.t. the distribution induced by the dynamics $\theta$ and the policies $\pi=(\pi_1,\ldots,\pi_K)$.

Similar to the previous section, we define a semi-distance on $\Theta$: \begin{equation}
    \label{dist2}
    d_\mu^2(\theta,\theta')=\int h^2(p_\theta(\cdot|s,\mathbf{a}),p_{\theta'}(\cdot|s,\mathbf{a}))\mu(s,\mathbf{a}) d\lambda(s,\mathbf{a}),
\end{equation}
where $\mu$ is a non-negative function and $\lambda$ is a reference finite measure on $\mathcal{S}\times\mathcal{A}_1\times\ldots\times\mathcal{A}_K$. The conditions for the existence of tests $\phi_t$ in a stochastic game is very similar to the single agent MDP. We can make slight modification to Assumption \ref{A4}, Assumption \ref{A5} and Theorem \ref{weakA4} by replacing $A_t$, $R_t$ with $\boldsymbol{A_t}$, $\boldsymbol{R_t}$ respectively, to find tests $\phi_t$ satisfying (\ref{typeI}) and (\ref{typeII}).  Then Theorem \ref{MDP} regarding the posterior convergence rate of $\theta$ still holds under the multi-agent scenario. Moreover, the posterior convergence rate for an episodic task in a stochastic game can also be obtained applying the same technique in the previous section and we omit the repeated details for simplicity. Furthermore, under Assumption \ref{A6}, we can easily obtain the convergence of the estimated V-functions for a fixed policy (\ref{Vpi}) by replacing $|V_\theta(\tilde{\pi})-V_{\theta_0}(\tilde{\pi})|$ with $|V_{k\theta}(\tilde{\pi}_1,\ldots,\tilde{\pi}_K)-V_{k\theta_0}(\tilde{\pi}_1,\ldots,\tilde{\pi}_K)|$, for $\forall\tilde{\pi}_1,\ldots,\tilde{\pi}_K\in \Pi$, $k=1,\ldots,K$.

According to the type of tasks, a stochastic game can be divided into three classes: fully cooperative, fully competitive and mixed \citep{busoniu2010reinforcement}. If $r_1=\ldots=r_K$ which means all agents aim to maximize the same expected return, the game is fully cooperative. If $K=2$ and $r_1=-r_2$, the game is fully competitive. The game is mixed if it is neither fully cooperative nor fully competitive. There exist many algorithms to find the optimal policies in a fully cooperative game including Team Q-learning \citep{littman2001value}, Distributed Q-learning \citep{lauer2000algorithm} and et al. In a fully cooperative game, if we obtain the optimal policies $\pi_{1\theta}^*,\ldots,\pi_{K\theta}^*$ under the dynamic $\theta$, we can obtain the convergence result (\ref{Vpiopt1}) by replacing $|V_{\theta}(\pi_\theta^*)-V_{\theta_0}(\pi_{\theta_0}^*)|$ with $|V_{k\theta}(\pi_{1\theta}^*,\ldots,\pi_{K\theta}^*)-V_{k\theta_0}(\pi_{1\theta_0}^*,\ldots,\pi_{K\theta_0}^*)|$ and the result (\ref{Vpiopt}) by replacing $|V_{\theta_0}(\pi_\theta^*)-V_{\theta_0}(\pi_{\theta_0}^*)|$ with $|V_{k\theta_0}(\pi_{1\theta}^*,\ldots,\pi_{K\theta}^*)-V_{k\theta_0}(\pi_{1\theta_0}^*,\ldots,\pi_{K\theta_0}^*)|$, $k=1,\ldots,K$. We should note that the convergence result of the regret bound (\ref{Vpiopt}) makes sense only if all the $K$ agents follow the optimal policies $\pi_{1\theta}^*,\ldots,\pi_{K\theta}^*$ computed under the same dynamics $\theta$ which means the $K$ agents should share the information of $\theta$ sampled from the posterior distribution $P_t(\cdot|H_t)$.

In a mixed stochastic game, an ordinary task for each agent is to find the Nash equilibrium (given in Definition \ref{Nash}) and apply the corresponding policy, denoted as $\pi_{k\theta}^*$ where $\theta$ is the dynamics of the game (Hu 1998). In the Nash equilibrium, each agent's policy is the best response to the other agents' policies. An adversarial equilibrium is a Nash equilibrium in which all agents are conflict with each other in some sense \citep{littman2001value}. It means an agent will achieve less return while the other agents will achieve more if the agent deviates from the equilibrium. The definition of an adversarial equilibrium is given in Definition \ref{adversarial}.
\begin{defn}
\label{Nash}
A Nash equilibrium in a stochastic game under the dynamics $\theta$ is a tuple of policies $(\pi_{1\theta}^*,\ldots,\pi_{K\theta}^*)$ such that 
\begin{equation}
\label{NashEqu}
\begin{aligned}
&V_{k\theta}(\pi_{1\theta}^*,\ldots,\pi_{K\theta}^*)\geq\\
&V_{k\theta}(\pi_{1\theta}^*,\ldots,\pi_{(k-1)\theta}^*,\pi_k,\pi_{(k+1)\theta}^*,\ldots,\pi_{K\theta}^*)
\end{aligned}
\end{equation}
for $\forall \pi_k\in\Pi$, $k=1,\ldots,K$.
\end{defn}
\begin{defn}
\label{adversarial}
An adversarial Nash equilibrium in a stochastic game under the dynamics $\theta$ is a tuple of policies $(\pi_{1\theta}^*,\ldots,\pi_{K\theta}^*)$ satisfying (\ref{NashEqu}) and
\begin{equation}
    \label{aNashEqu}
    V_{k\theta}(\pi_{1\theta}^*,\ldots,\pi_{K\theta}^*)\leq V_{k\theta}(\pi_1,\ldots,\pi_{k-1},\pi_{k\theta}^*,\pi_{k+1},\ldots,\pi_K)
\end{equation}
for $\forall \pi_i\in \Pi$, $i\neq k$, $k=1,\ldots,K$.
\end{defn}

In a mixed stochastic game, especially an adversarial one, the agents may not share their information including their current policies and their estimates for the dynamics $\theta$.
In order to learn the convergence rates of the estimated V-functions and the regret bound in a mixed game, we make Assumption \ref{A7} considering the structure of the game.
\begin{assumption}
\label{A7}
Suppose $1\leq k\neq j\leq K$ and $\pi_k\neq \pi_{k\theta}^*$. Then we have
$$
V_{k\theta}(\pi_{1\theta}^*,\ldots,\pi_{K\theta}^*)\geq V_{k\theta}(\pi_1,\ldots,\pi_k,\ldots,\pi_{j\theta}^*,\ldots,\pi_K)
$$
for $\forall \pi_i\in \Pi$, $i\neq k,j$.
\end{assumption}
Assumption \ref{A7} means that if an agent applies his Nash equilibrium policy, the best response for another agent is also his Nash equilibrium policy, no matter what the other agents' policies are. Notice that the assumption is automatically satisfied if $K=2$. 

We consider the convergence rate of the regret bound for an agent's V-function in two cases: (i) the other agents are following their Nash equilibrium policies (ii) for all $k=1,\ldots,K$, the $k$-th agent applies the estimated optimal policy $\pi_{k\hat{\theta}_k}^*$ where $\hat{\theta}_k$ is sampled from $P_t(\cdot|H_t)$. The convergence results are shown in Corollary \ref{multV}.
\begin{corollary}
\label{multV}
If Assumption \ref{A6} and (\ref{conv}) hold, then for every $M_t\rightarrow \infty$, we have
\begin{equation}
\begin{aligned}
&P_t(\theta\in \Theta: |V_{k\theta_0}(\pi_{1\theta_0}^*,\ldots,\pi_{k\theta}^*,\ldots,\pi_{K\theta_0}^*)-\\
&V_{k\theta_0}(\pi_{1\theta_0}^*,\ldots,\pi_{K\theta_0}^*)|> M_t\epsilon_t^\rho|H_t)\rightarrow 0
\end{aligned}
\end{equation}
in $P_{\theta_0,\pi}$-probability for $\forall k=1,\ldots,K$.\\
Moreover, if Assumption \ref{A7} holds, then we have
\begin{equation}
\begin{aligned}
&P_{\theta_0,\pi}(\underset{k=1,\ldots,K}{\max}|V_{k\theta_0}(\pi_{1\hat{\theta}_1}^*,\ldots,\pi_{K\hat{\theta}_K}^*)-\\
&V_{k\theta_0}(\pi_{1\theta_0}^*,\ldots,\pi_{K\theta_0}^*)|> M_t\epsilon_t^\rho)\rightarrow 0
\end{aligned}
\end{equation}
as $t\rightarrow \infty$.
\end{corollary}

\section{Numerical Results}
\subsection{Toy Example}
We consider a toy example with $\mathcal{S}=\{0,1\}$, $\mathcal{A}=\{0,1\}$. The transition dynamics and reward function are given as follows:
\begin{enumerate}[(i)]
    \item If $S_t=1,A_t=1$, then $S_{t+1}\sim \text{Bin}(1,\theta_1)$, $r(S_t,A_t)=1$.
    \item If $S_t=1,A_t=0$, then $S_{t+1}\sim \text{Bin}(1,1-\theta_1)$, $r(S_t,A_t)=0.5$.
    \item If $S_t=0,A_t=1$, then $S_{t+1}\sim \text{Bin}(1,\theta_2)$, $r(S_t,A_t)=1.5$.
    \item If $S_t=0,A_t=0$, then $S_{t+1}\sim \text{Bin}(1,1-\theta_2)$, $r(S_t,A_t)=2$.
\end{enumerate}
We set the true parameter $\theta_{10}=0.2$, $\theta_{20}=0.4$ and $\Theta=\{(\theta_1,\theta_2):0.01\leq \theta_1,\theta_2\leq 0.99\}$. The priors for $\theta_1$ and $\theta_2$ are $U(0.01,0.99)$ and they are independent. We set the length of MDPs $T=5000$, the initial state $s_0=0$ and the discounted factor $\gamma=0.25$. We study the posterior convergence rate of the $(\theta_1,\theta_2)$, the estimated V-function and the regret bound under two scenarios: (i) the agent select the actions $0,1$ with equal probabilities at each time step; (ii) the agent applies Algorithm \ref{TS} with $\delta_t=t^{-\frac{1}{4}}$. According to Theorem \ref{weakA4} and Theorem \ref{MDP}, we can obtain the posterior convergence rates $\epsilon_t$ in (\ref{conv}) for the above two scenarios are $t^{-\frac{1}{2}+\delta}$ and $t^{-\frac{1}{4}+\delta}$ respectively, for $\forall \delta>0$.  Moreover, we can show the convergence rates of $(\theta_1,\theta_2)$ w.r.t. $L_2$-norm, the estimated V-function and the regret bound for the two scenarios are also $t^{-\frac{1}{2}+\delta}$ and $t^{-\frac{1}{4}+\delta}$ respectively (technical details in Appendix). Figure \ref{fig:toy} displays the average difference between the posterior samples $(\theta_1,\theta_2)$ and $(\theta_{10},\theta_{20})$ w.r.t. $L_2$-norm, the average difference between the estimated optimal V-function and the true one and the average regret bound at each time step. $100$ Monte Carlo runs are simulated. In Figure \ref{fig:toy}, we see that the theoretical convergence rates given above are verified and the three quantities usually converges faster than the theoretical results.

\subsection{RiverSwim Example}
\begin{figure}
  \centering
  \includegraphics[width=0.5\textwidth]{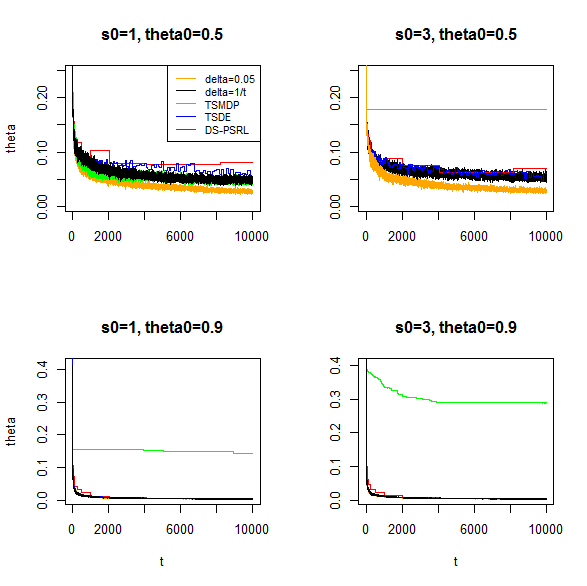}
  \caption{\label{fig:swim} The average difference between posterior samples $\theta$ and $\theta_0$ for $\epsilon$-greedy Thompson sampling with $\delta_t=0.05$, $\epsilon$-greedy Thompson sampling with $\delta_t=1/t$, TSMDP, TSDE and DS-PSRL.
  }
\end{figure}

We consider the RiverSwim example \citep{strehl2008analysis} which models an agent swimming in a river who can choose to swim either left or right. The MDP consists of six states arranged in a chain with the agent starting in the state $s_0$. If the agent decides to move left i.e. with the river current, he is always successful; if he decides to move right, he might succeed with probability $\theta$ otherwise he stays in the current state. The reward function is given by: $r(s,a)=2$ if $s=1$ and $a=\text{left}$; $r(s,a)=10$ if $s=6$ and $a=\text{right}$; otherwise $r(s,a)=0$. We set $\delta_t=t^{\frac{1}{2}}$, $\delta_t=0.05$ in Algorithm \ref{TS} and compare our method to three Thompson sampling algorithms: TSMDP, TSDE and DS-PSRL. We let $\Theta=\{\theta:0.01\leq \theta\leq 0.99\}$ and the prior $U(0.01,0.99)$. Let the length of MDPs $T=10000$ and the discounted factor $\gamma=0.99$. We consider the proportion of selecting the optimal action and the convergence of the dynamics $\theta$ as the measures of performance. Simulations are conducted in four cases: (i) $s_0=3$, $\theta_0=0.5$ (ii) $s_0=3$, $\theta_0=0.9$ (iii) $s_0=1$, $\theta_0=0.5$ (iv) $s_0=1$, $\theta_0=0.9$ and $100$ Monte Carlo runs are simulated. Table \ref{tab:performance_comparison} shows the average proportion of selecting the optimal action for the five methods under the four scenarios. We can see that $\epsilon$-greedy Thompson sampling with $\delta_t=1/t$ always performs better than DS-PSRL. In the case $s_0=3$ and $\theta_0=0.9$, TSMDP achieves a very low proportion of selecting the optimal action since $s_0=3$ is not recurrent under some sub-optimal policy. Thus the performance of TSMDP depends heavily on an appropriate choice of the initial state which is sometimes not possible for a general unknown MDP. TSDE can always achieve a high proportion of selecting the optimal action in this example but the method is limited to finite MDPs. Figure \ref{fig:swim} displays the average difference between the posterior samples $\theta$ and $\theta_0$ at each time step for the five methods. We can see that the posterior samples $\theta$ always converge faster to $\theta_0$ in $\epsilon$-greedy Thompson sampling algorithms with $\delta=0.05$ and $\delta=1/t$  than in other methods which give no theoretical guarantee. Combining with Table \ref{tab:performance_comparison}, we can see the trade-off between selecting the optimal action and learning the dynamics $\theta$.
\subsection{Glucose Example}
In this experiment, we consider an MDP in which the states are continuous. We simulate cohorts of patients with type 1 diabetes using a generative model based on the mobile health study of \cite{maahs2012outpatient}. We only consider the action of whether to use insulin, so the action space $\mathcal{A}=\{0,1\}.$ The covariates observed for patient $i$ at time $t$ is average blood glucose level, total dietary intake, and total counts of physical activity, denoted by $(Gl_{i,t},Di_{i,t},Ex_{i,t})$ respectively. Glucose levels evolve according to the AR(2) process
$$
\begin{aligned}
Gl_t&=\beta_0+\beta_1 Gl_{t-1}+\beta_2 Di_{t-1}+\beta_3 Ex_{t-1}+\beta_4 Gl_{t-2}\\
&+\beta_5 Di_{t-2}+\beta_6 Ex_{t-2}+\beta_7 A_{t-2}+\beta_8 A_{t-1}+e_t,
\end{aligned}
$$
where $e_t\sim N(0,\sigma^2)$, $\sigma=5$ and $(\beta_0,\beta_1,\ldots,\beta_8)=(10,0.9,0.1,-0.01,0.0,0.1,-0.01,-10,-4)$; $A_t$ is the action taken at time $t$; and $Di_t\sim N(\mu_d,\sigma^2_d)$ with probability $p_d$, otherwise $Di_t=0$; similarly, $Ex_t\sim N(\mu_e,\sigma_e^2)$ with probability $p_e$, otherwise $Ex_t=0$, where $\mu_d=\mu_e=0$, $\sigma_d=\sigma_e=10$, $p_d=p_e=0.6$. Thus the dynamics are Markovian with states $S_t=(Gl_t,Di_t,Ex_t,Gl_{t-1},Di_{t-1},Ex_{t-1},A_{t-1})$. The reward at each time step is given by $R_t=1(Gl_t<70)[-0.005Gl_t^2+0.95Gl_t-45]+1(Gl_t\geq 70)[-0.0002Gl_t^2+0.022Gl_t-0.5]$, which decreases as $Gl_t$ departs from normal glucose levels.

In our experiments, we simulate data for 70 patients and time horizons $T=30$ and $T=50$. We assume $\boldsymbol{\beta}=(\beta_0,\ldots,\beta_8)'$ are the only unknown parameters and the prior for $\boldsymbol{\beta}$ is $N(\mathbf{0},\frac{1}{4}\mathbf{I}_9)$. We apply Algorithm \ref{TS} with $\delta_t=0.05$ and consider the cumulative rewards up to $T$ as the measure of performance, averaged over 50 Monte Carlo runs. For a posterior sample $\boldsymbol{\beta}$ at each time step, we compute the optimal policy with the algorithm of Fitted Q Iteration (FQI) which is a tree-based method \citep{ernst2005tree}. For the FQI algorithm, we generate sufficient training sets, i.e. sets of four tuples $\{(s_t,a_t,r_t,s_{t+1})\}_{t\geq 0}$, under the sampled dynamics $\boldsymbol{\beta}$, set the number of iterations $N=5$ and fit with random forest regression. We compare Algorithm \ref{TS} with three other methods: DS-PSRL, the gold standard method and naive FQI. The gold standard method assumes $\beta_0,\ldots,\beta_8$ are completely known and compute the optimal policy with FQI by generating sufficient training sets under the true dynamics. The naive FQI is a model-free method which only uses the past data generated up to current time step, i.e. the history $H_t$, as the training sets to fit the Q-function. Table \ref{tab:glucose} shows the average cumulative rewards when $T=30,50$ for the above four methods. We can see that $\epsilon$-greedy Thompson sampling performs better than DS-PSRL and the naive FQI, which results from the fact that our method updates the dynamics $\beta$ more frequently than DS-PSRL and employs knowledge of the model structure while the naive FQI is model-free. As we expect, the gold standard method always performs the best.
\section{Discussion}
In this paper, we show the convergence rates of posterior distributions of the model dynamics in single-agent MDPs and Markov games for both episodic and continuous tasks. The theoretical results hold for general state and action space and the parameter space of the dynamics can be infinite dimensional. Moreover, we show the convergence rates of posterior distributions of the mean accumulative reward under a fixed or the optimal policy and of the regret bound. Then we propose the $\epsilon$-greedy Thompson sampling algorithm which provides both posterior convergence rates for the dynamics and the regret-type bound. The numerical results verify the validity of convergence rates the theorems give and the competitiveness of our proposed Thompson sampling algorithm compared to others. 

We only consider the transition density as the model dynamics in this paper but it will not need additional techniques to include the reward function. For the future work, we can explore conditions under which our convergence results still hold for the completely greedy Thompson sampling algorithm. Moreover, we can extend our current results to a partially observable MDP.

\begin{figure}
  \centering
  \includegraphics[width=0.5\textwidth]{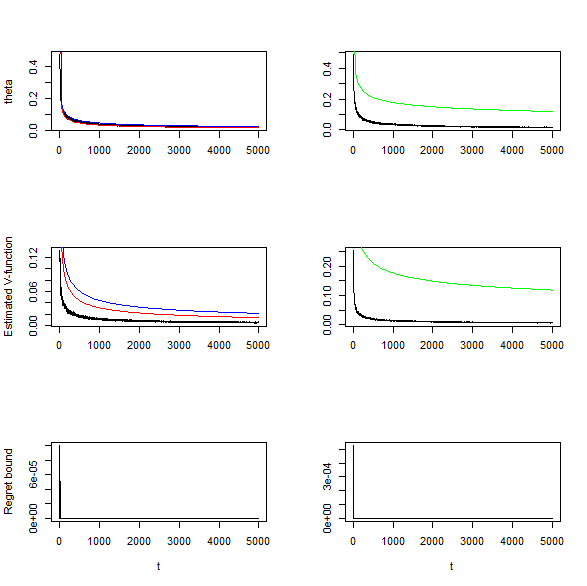}
  \caption{\label{fig:toy} The results of convergence rates in scenarios (i) and (ii). The left and right column display the results for scenarios (i) and (ii) respectively. The black lines in the three rows display the average difference between posterior samples $(\theta_1,\theta_2)$ and $(\theta_{10},\theta_{20})$, the average difference between the estimated optimal V-function and the true one and the average regret bound respectively. The red, blue and green curve represent $f(t)=t^{-\frac{1}{2}}$, $f(t)=t^{-\frac{9}{20}}$ and $f(t)=t^{-\frac{1}{4}}$ respectively.
  }
\end{figure}
\begin{table}[tp]  
  \centering  
  \fontsize{6.5}{8}\selectfont  
  \begin{threeparttable}  
  \caption{The average proportion of selecting the optimal action for five methods under different scenarios. The standard errors are given in the brackets.}  
  \label{tab:performance_comparison}  
    \begin{tabular}{ccccc}  
    \toprule  
    \multirow{2}{*}{Method}&  
    \multicolumn{2}{c}{ $s=1$}&\multicolumn{2}{c}{$s=3$}\cr  
    \cmidrule(lr){2-3} \cmidrule(lr){4-5}  
    &$\theta_0=0.5$&$\theta_0=0.9$&$\theta_0=0.5$&$\theta_0=0.9$\cr  
    \midrule  
    $\epsilon$-greedy TS ($\delta_t=.05$)&$.9628 (10^{-5})$&$.9747 (10^{-5})$&$.9622 (10^{-4})$&$.9748 (10^{-5})$\cr 
    $\epsilon$-greedy TS ($\delta_t=1/t$)&$.9836 (10^{-4})$&$.9993 (10^{-6})$&$.9846 (10^{-4})$&$.9993 (10^{-6})$\cr  
    TSMDP&$.9872 (10^{-4})$&$.9998 (10^{-6})$&$.9998 (10^{-5})$&$.57 (.005)$\cr  
    TSDE&$.9897 (10^{-4})$&$.9994 (10^{-5})$&$.9894 (10^{-4})$&$.9995 (10^{-5})$\cr  
    DS-PSRL&$.9740 (10^{-4})$&$.9976 (10^{-4})$&$.9738 (10^{-4})$&$.9980 (10^{-4})$\cr  
    \bottomrule  
    \end{tabular}  
    \end{threeparttable}  
\end{table}
\begin{table}[!htbp]  
  \centering  
  \fontsize{6.5}{8}\selectfont  
  \caption{The average cumulative reward when $T=30,50$ for four methods. The standard errors are given in the brackets.}  
  \label{tab:glucose}  
    \begin{tabular}{ccc}
    \hline
    Method &$T=30$ & $T=50$\\
    \hline
    $\epsilon$-greedy TS ($\delta_t=.05$)&-27.68 (0.55)&-38.91 (0.74)\\
    DS-PSRL&-28.17 (0.62)&-46.05 (1.30)\\
    Gold standard &-18.48 (0.20)&-23.18 (0.35) \\
    Naive FQI &-34.72 (1.24)&-41.68 (1.33)\\
    \hline
    \end{tabular}  
\end{table}
\clearpage
\bibliographystyle{aaai}
\bibliography{mybib}

\end{document}